\documentclass[11pt,a4paper]{article} 
\usepackage{amsmath, amssymb,amscd, amsthm} 
\title{Integrality and the Laurent phenomenon for Somos $4$ sequences} 
\author{Christine Swart\thanks{Department of Mathematics \& 
Applied Mathematics, University 
of Cape Town, Rondebosch, 7700, South Africa. ~~E-mail: cswart@maths.uct.ac.za} 
$\,$ and Andrew Hone\thanks{Institute of Mathematics, 
Statistics and Actuarial Science, University of
Kent, Canterbury CT2 7NF, U.K. ~~E-mail: A.N.W.Hone@kent.ac.uk}
} 
 
\begin{document} 
 
\renewcommand{\theequation}{\arabic{section}.\arabic{equation}}

\newcommand{\beq}{\begin{equation}}  
\newcommand{\eeq}{\end{equation}}  
\newcommand{\bea}{\begin{eqnarray}}  
\newcommand{\eea}{\end{eqnarray}}  
\newcommand\la{{\lambda}}   
\newcommand\ka{{\kappa}}   
\newcommand\al{{\alpha}}   
\newcommand\be{{\beta}}   
\newcommand\om{{\omega}}  
\newcommand\tal{{\tilde{\alpha}}}  
\newcommand\tbe{{\tilde{\beta}}}   
\newcommand\tla{{\tilde{\lambda}}}  
\newcommand\tmu{{\tilde{\mu}}}  
\newcommand\si{{\sigma}}  
\newcommand\lax{{\bf L}}    
\newcommand\mma{{\bf M}}    
\newcommand\rd{{\mathrm{d}}}  
\newcommand\tJ{{\tilde{J}}}  
\newcommand\I{{\mathcal{I}}}
\newcommand\tI{{\tilde{\mathcal{I}}}}

\newcommand\tis{{\tilde{s}}} 
\newcommand\tit{{\tilde{t}}}
\newcommand\tiu{{\tilde{u}}}
\newcommand\tiv{{\tilde{v}}}
\newcommand\tiw{{\tilde{w}}}

\newtheorem{propn}{Proposition}[section]
 
\newtheorem{thm}{Theorem}[section]

\newtheorem{exa}[thm]{Example}

\newtheorem{lem}{Lemma}[section]

\newtheorem{cor}[thm]{Corollary}
\newenvironment{prf}{\trivlist \item [\hskip
\labelsep {\bf Proof.}]\ignorespaces}{\qed \endtrivlist}

\newenvironment{rem}{\trivlist \item [\hskip

\labelsep {\bf Remark:}]\ignorespaces}{\endtrivlist}

\newcommand{\N}{{\mathbb N}}

\newcommand{\Q}{{\mathbb Q}}

\newcommand{\Z}{{\mathbb Z}}

\newcommand{\C}{{\mathbb C}}

\maketitle

\begin{abstract} 
Somos 4 sequences  
are  a family of sequences defined by a fourth-order 
quadratic recurrence relation with constant coefficients. 
For particular choices of the coefficients and 
the four initial data, such recurrences can  
yield sequences of integers. Fomin and Zelevinsky 
have used the theory of cluster 
algebras to prove that these recurrences also provide one 
of the simplest examples of the Laurent phenomenon: 
all the terms of a Somos 4 sequence are Laurent 
polynomials in the initial data. The integrality 
of certain Somos 4 sequences has previously been understood in terms of  
the Laurent phenomenon. However, each of the 
authors of this paper has independently established 
the precise correspondence between 
Somos 4 sequences and sequences of points on elliptic curves.    
Here we show that these sequences 
satisfy a  stronger 
condition than the Laurent property, and hence establish 
a broad set of sufficient conditions for integrality. 
As a by-product, non-periodic 
sequences provide infinitely many solutions of an associated 
quartic Diophantine equation in four variables.
The  analogous 
results for Somos 5 sequences are also presented, 
as well as various  
examples, including parameter families of Somos 4 integer sequences.  
\end{abstract} 

\section{Introduction}

The study of integer sequences generated by  
linear recurrences has a long history, 
and, apart from their intrinsic interest in number theory, 
nowadays there are many applications of such sequences 
in computer science and cryptography. However, 
not so much is known about sequences generated by nonlinear recurrences. 
In  \cite{recs} (see e.g. section 1.1.20) it is suggested that 
Somos sequences, which are generated by quadratic recurrence 
relations of the form  
$$
S_{n+k}S_n=\sum_{j=1}^{[k/2]}\al_j \, S_{n+k-j}S_{n+j}
$$ 
(where the coefficients $\al_j$ are constant), are 
suitable generalizations of linear recurrence sequences, in the sense 
that they have many analogous properties. In this paper we are mainly 
concerned with the cases $k=4$ and $k=5$, and the problem of determining 
when such recurrences yield sequences of integers; unlike 
the  situation for linear recurrences, this problem is not so straightforward. 

Somos $4$ sequences are defined by a fourth-order quadratic 
recurrence relation of the form 
\beq 
A_{n+4}
A_{n}
=\al \,
A_{n+3}  
A_{n+1}  
+\beta  
\,  
(A_{n+2})^2,  
\label{bil}  
\eeq  
where $\al$, $\be$ are constant coefficients. We also 
refer to equations like (\ref{bil}) as {\it bilinear} 
recurrences, by analogy with the Hirota bilinear form 
of integrable partial differential equations in soliton 
theory \cite{hir}, which have discrete (difference 
equation) counterparts  \cite{zabrodin}.  
During an exploration of the combinatorial properties 
of elliptic theta functions, 
Michael Somos introduced the sequence defined by (\ref{bil}) 
with coefficients $\al = \be =1$ and initial data 
$A_1=A_2=A_3=A_4=1$, which begins 
\beq \label{s4seq} 
1,1,1,1,2,3,7,23,59,314,1529, 8209, 83313,\ldots 
\eeq 
(and in this case the sequence extends symmetrically backwards 
as $\ldots ,59,23,$ 
$7,3,2,1,1,1,1,2,\ldots$). It is a remarkable 
fact that although the recurrence defining the 
sequence (\ref{s4seq}) is rational 
(because one must divide by $A_n$ in order to obtain $A_{n+4}$), 
this sequence consists entirely of integers.  
In what follows we are  
mainly concerned with Somos 4 sequences with integer 
coefficients and initial data, and we  establish a broad set 
of sufficient conditions for integrality of such sequences. 
However, it is most convenient 
to work with sequences over $\C$ and indicate when   
our results reduce to sequences over $\Q$ or to integer sequences. 
Given four adjacent non-zero values, the recurrence (\ref{bil}) 
can be iterated either forwards or backwards, so it is natural 
to define the sequences for $n\in\Z$. 
(In fact, if a zero is encountered then 
the sequence can be analytically continued through it, either by using the 
fact that the recurrence (\ref{bil}) has the singularity 
confinement property of \cite{grp}, or by making use of the 
explicit formula for the iterates given in (\ref{form}) below.) 

The particular sequence (\ref{s4seq}) above is sometimes referred to 
as {\it the} Somos 4 sequence \cite{sloane}, and the first published 
proof of its integrality appears in the survey by
Gale \cite{gale}.  
However, following \cite{rob} we will refer to 
(\ref{s4seq}) as the Somos (4) sequence, in order to distinguish it 
from other 
Somos 4 sequences defined by 
recurrences of the general form (\ref{bil}). Such sequences originally arose in 
the theory of elliptic 
divisibility sequences (EDS), which were introduced by Morgan Ward 
\cite{ward1, ward2}, who considered integer sequences satisfying 
a recurrence of the form 
\beq \label{eds} 
W_{n+4}W_n=(W_2)^2 \, W_{n+3}\,W_{n+1} -W_1\,W_3\,(W_{n+2})^2 , 
\eeq 
these being a special case of (\ref{bil}) with 
\beq \label{edsp} 
\al = (W_2)^2, \qquad \be = -W_1\,W_3. 
\eeq 
The form of the recurrence (\ref{eds}) defines 
an antisymmetric sequence, so that $W_n=-W_{-n}$. 
Ward showed \cite{ward1} that 
if the initial data are  taken to be 
\beq \label{edsinits} 
W_0=0, \quad W_1=1, \quad W_2,W_3,W_4\in\mathbb{Z}  
\qquad \mathrm{with}\quad W_2|W_4,  
\eeq   
then it turns out that the subsequent terms of the sequence 
are all integers satisfying the 
divisibility property
\beq  
\label{divis}  
W_n  |W_m \qquad \mathrm{ whenever}  
\qquad n|m.  
\eeq 
In this sense, EDS constitute a nonlinear generalization of 
Lucas sequences. 
Furthermore, Ward demonstrated 
that the terms of a generic EDS defined by (\ref{eds}) 
correspond to the 
multiples $[n]P$ of a point $P$ on an elliptic curve $E$, and the 
$n$th term in the sequence can be written explicitly 
in terms of the Weierstrass sigma function associated with the curve $E$, 
as 
\beq 
W_n=\frac{\si (n\ka)}{\si (\ka )^{n^2}}.  
\label{wardeds} 
\eeq 

Somos 4 sequences, and EDS in particular, are of considerable 
interest to number theorists due to the way that 
primitive prime divisors appear  
\cite{eew, ems}. In that and in 
other respects, they have a lot in common with linear recurrence 
sequences (see \cite{recs}, chapters 10 \& 11). 
Somos 4 sequences have also arisen 
quite recently 
in a different guise, providing  one of the simplest examples 
of the combinatorial structures appearing in Fomin and Zelevinsky's 
theory of cluster algebras. More precisely, the following 
result was proved in \cite{fz}. 

\begin{thm} {\bf The Laurent property for Somos 4 sequences:} \label{lprops4}
For a Somos 4 sequence defined by a fourth-order recurrence 
(\ref{bil}) with coefficients $\al$, $\be$ 
and initial data $A_1, A_2, A_3, A_4$, all of the terms in the sequence 
are Laurent polynomials in these initial data whose coefficients 
are in $\Z [\al ,\be ]$, so that $A_n \in \Z [\al , \be , A_1^{\pm 1}, 
A_2^{\pm 1},  A_3^{\pm 1},  A_4^{\pm 1}]$ for all $n\in \Z$.         
\end{thm}  

The above result is only one of a series of analogous 
statements  for a wide variety 
of recurrences in one and more dimensions to which 
the machinery of cluster algebras was applied in \cite{fz}. The integrality 
of the  
Somos (4) sequence (\ref{s4seq}), and also of any other 
sequence defined by the recurrence (\ref{bil}) with coefficients 
$\al$, $\be\in\Z$ and initial values $A_1=A_2=A_3=A_4=1$, 
is an immediate consequence of the Laurent property.   However, 
there are many other examples of Somos 4 sequences that 
consist entirely of integers, but for which the 
Laurent phenomenon alone is insufficient to deduce integrality. 
For instance, consider the Somos 4 sequence generated by 
\beq\label{1stex}
A_{n+4}A_{n} = 1331  A_{n+3}A_{n+1} + 119790 A_{n+2}^2  
\eeq 
with the initial values 
$$A_1=1, \quad A_2=3, \quad A_3=121, \quad 
A_4=177023. 
$$ 
This turns out to be an integer sequence, i.e. $A_n\in\Z$ for all $n\in\Z$, but Theorem
\ref{lprops4} is not enough to show this; it follows from Theorem \ref{strongs4} below, which 
relies on the connection between Somos  4 sequences and elliptic curves.

In the rest of the paper we 
shall make use of the fact that the terms $A_n$ in a 
Somos 4 sequence correspond to 
a sequence of points $Q+[n]P$ on an associated elliptic curve $E$, 
which was  proved by 
each  of us independently, in \cite{swart}
and \cite{honeblms} respectively. Our 
results can be combined and summarized as follows. 

\begin{thm} \label{swahonresults} 
If $(A_n)$ is a Somos 4 sequence defined by a  recurrence 
(\ref{bil}) with coefficients $\al$, $\be$, then the 
quantity 



\beq 
\label{tinvt} 
T = \dfrac{A_{n-1} \, A_{n+2}}{A_n \, A_{n+1}} + \al \left( 
\dfrac{ (A_n)^2}{A_{n-1} \, A_{n+1}} + \dfrac{(A_{n+1})^2}{A_n\, A_{n+2}} 
\right)  
    + \dfrac{\be\,A_n\,A_{n+1} }{A_{n-1} \, A_{n+2}} 
\eeq  
is independent of $n$, and for $\al \neq 0$ the 
sequence corresponds to a sequence of points $Q+[n]P$ on an 
elliptic curve $E$ with $j$-invariant 
\beq 
\label{jinvt} 
j=\frac{(T^4-8\be T^2  -24\al^2 T +16\be^2)^3}{\al^4 (\be T^4+\al^2T^3 
-8\be^2T^2-36\al^2\be T+16\be^3-27\al^4)}  
\eeq 
Moreover, the general term of the sequence can be written 
in the form 
\beq
A_n=
\hat{a}\,\hat{b}^n\,\frac{\si (z_0 +n\ka )}{
\si (\ka )^{n^2}},
\label{form}
\eeq
where $\si (z)=\si (z;g_2,g_3)$  
denotes the Weierstrass sigma function associated with the 
curve $E$ written in the canonical form  
\beq \label{wcurve} 
E: \qquad y^2=4x^3-g_2x-g_3, 
\eeq 
and the six parameters $g_2,g_3,\hat{a},\hat{b},z_0,\ka \in\C$ are 
determined uniquely (up to a choice of sign) 
from the two coefficients 
$\al ,\be $ and the four initial 
data $A_1,A_2,A_3,A_4$.  The $(x,y)$ coordinates of 
the corresponding sequence of points on the curve $E$ 
are given by  
$Q+[n]P=(\,\wp (z_0+n\ka )\, , \,\wp '(z_0+n\ka )\,)$.   \end{thm} 

\noindent {\bf Remarks.} 
Given four non-zero initial data $A_1,A_2,A_3,A_4$, 
the quantity $T$ in (\ref{tinvt}) can be calculated 
as 
\beq \label{tinvt2}  
T = \dfrac{(A_1)^2 \, (A_4)^2 + \alpha ((A_2)^3 \, A_4 + A_1 \, (A_3)^3)  
+ \beta \, (A_2)^2 \, (A_3)^2}{A_1 \, A_2 \, A_3 \,
A_4}. \eeq  
(For a discussion of zeros in Somos sequences, see \cite{swart}.) 
The translation invariant was first 
found in \cite{swart}, and also appeared in \cite{hones5}, where it was 
denoted $J$,  and written in terms of  
$ 
f_n=A_{n-1}\,A_{n+1} / (A_n)^2  
$ 
as 
\beq \label{tmap} 
T=f_n\,f_{n+1}+\al \left( \frac{1}{f_n}+\frac{1}{f_{n+1}}\right) 
+\frac{\be}{f_n\,f_{n+1}}. 
\eeq   
(Note that 
the equation $T=\mathrm{constant}$ defines a quartic curve 
in the $(f_n,f_{n+1})$  
plane: this curve is birationally equivalent to the curve 
$E$ given in (\ref{wcurve}).)
In the paper \cite{honeblms}, a more complicated 
invariant was used, namely the quantity  
\beq 
\label{la} 
\la =(T^2/4-\be )/(3\al ) . 
\eeq     
We can also use the equation (\ref{bil}) to  write 
$T$ as a more compact expression in terms of five adjacent terms 
of the sequence, that is 
\beq
\label{tinvt3}
T = \dfrac{A_{n-1} \, A_{n+2}}{A_n \, A_{n+1}} + \dfrac{\alpha \, A_n^2}{A_{n-1}
 \, A_{n+1}}  
    + \dfrac{A_{n-2} \, A_{n+1}}{A_{n-1} \, A_n}
\eeq
The explicit formulae for the invariants $g_2$, $g_3$ of the curve 
(\ref{wcurve}) in terms of $\al$, $\be $ and $T$  
are given in section 2 below, along with a summary 
of other results that are needed here. 
It is worth mentioning that the above theorem  
still holds in the degenerate case when the bracketed expression 
in the numerator of (\ref{jinvt}) vanishes, so that $j=\infty$   
and the curve $E$ becomes singular; in that case the 
formula (\ref{form}) is written in terms of 
hyperbolic functions (or, when 
$g_2=g_3=0$, in terms of rational functions).

The connection between Somos 4 and Somos 5 sequences 
and elliptic curves was previously understood in 
unpublished work of several number theorists: see Zagier's discussion of Somos 5 \cite{zagier}, 
the results of Elkies quoted in \cite{heron}, and the unpublished work of Nelson Stephens 
mentioned in \cite{swart},  
where
alternative Weierstrass 
models for the elliptic curve $E$ are used. 
For 
another approach using birationally equivalent 
quartic curves and continued fractions, see \cite{vdp}. 
A brief history of Somos sequences can be 
found on Propp's Somos sequence site \cite{propp}.

\section{Companion elliptic divisibility sequences} 

\setcounter{equation}{0}

The purpose of this section is to present some facts about 
Somos 4 sequences and their relationship with elliptic curves 
and division polynomials which will be useful for what follows. 
Given a Somos 4 
sequence with coefficients $\al$, $\be$ 
and initial data $A_1,A_2,A_3,A_4$ we can immediately compute 
the translation invariant $T$ given by (\ref{tinvt}), and then the 
invariants for the associated curve (\ref{wcurve}) can be 
calculated directly from the formulae 
\beq \label{g2} 
g_2= \frac{T^4-8\be T^2-24\al^2 T+16\be^2}{12\al^2},  
\eeq   
\beq \label{g3}
g_3=-\frac{T^6-12\be T^4-36\al^2 T^3+48\be^2T^2+144\al^2\be T + 
216\al^4-64\be^3}{216\al^3},  
\eeq
while $j=1728g_2^3/(g_2^3-27g_3^2)$ yields the expression (\ref{jinvt}) for the $j$-invariant. 
(Upon substituting for $\la$ 
as in (\ref{la}),  
the formulae (\ref{g2}) and (\ref{g3}) are seen to be 
equivalent to the 
formulae given in \cite{honeblms}.)

The backward  
iterates $A_0$, $A_{-1}$ and $A_{-2}$ are obtained by applying (\ref{bil}) 
in the reverse direction, so that the quantities 
$$ 
f_{-1} =A_{-2}\,A_{0}/(A_{-1})^2, \quad 
f_0=A_{-1}\,A_{1}/(A_0)^2,
\quad 
f_1=A_{0}\,A_{2}/(A_1)^2,
\qquad
$$ 
can be calculated, and then the associated 
sequence of points on the curve  
(\ref{wcurve}) is given by  
\beq 
\label{pts} \begin{array}{ccl}   
Q+[n]P & =  & \Big(\,\wp (z_0)\, , \, \wp '(z_0)\, \Big)+[n]\, 
\Big(\,\wp (\ka )\, , \, \wp '(\ka )\, \Big) \\ 
&& \\ 
& = &  
\Big(\,\la -f_0\, , \, (f_0)^2(f_1-f_{-1})/\sqrt{\al}\, \Big) 
+[n]\, \Big(\, \la \, , \, \sqrt{\al}\, \Big) \end{array}  
\eeq
with $\la$ as in  (\ref{la}). Up to the overall ambiguity in the sign 
of $\sqrt{\al}=\wp'(\ka )$, which corresponds to the freedom to send 
$\ka\to -\ka$, $z_0\to -z_0$ on $\mathrm{Jac}(E)$, the 
coefficients 
and the translation invariant are given as elliptic functions 
of $\ka$ by 
\beq 
\label{parform} 
\al =\wp '(\ka )^2, \quad \be =\wp '(\ka )^2\Big(\wp (2\ka )-\wp (\ka )\Big), 
\quad T=\wp ''(\ka ). 
\eeq

Elliptic divisibility sequences (EDS) 
were originally defined by Morgan Ward \cite{ward1} 
as 
sequences of integers satisfying the property (\ref{divis}) as 
well as the family of recurrences 
\beq \label{hankel}
W_{n+m}\, W_{n-m}=\left|\begin{array}{cc} W_m W_{n-1} &
W_{m-1}W_n \\
W_{m+1}W_{n} &
W_{m}W_{n+1} \end{array}\right|,
\eeq
for all $m,n$, and they were called {\it proper} EDS if 
$W_0=0$, $W_1=1$ and $W_2W_3\neq 0$. 
Here we define a {\it generalized} EDS to 
be a sequence $\{\,W_n \,\}_{n\in\Z }$ 
specified by four non-zero initial data 
$$ 
W_1=1, \quad W_2,W_3,W_4\in\C^* 
$$ 
together with the recurrence (\ref{eds}), 
of which immediate 
consequences are that $W_0=0$ and the sequence is consistently 
extended to negative indices so that $W_{-n}=-W_n$ and 
(\ref{eds}) holds for all $n\in\Z$. We may then make the 
observation that a   
(generalized) EDS  is just a special type of Somos 4 sequence, and thus
it follows immediately from Theorem 1.2 (i.e. the 
formula (\ref{form}) with $\hat{a}=\hat{b}=1$ and $z_0=0$) that the terms of 
the EDS are given by the explicit expression (\ref{wardeds}) 
in terms of the Weierstrass sigma function. 
The whole family of recurrences (\ref{hankel}) 
then follows as a direct consequence of this 
explicit formula together with the three-term equation for 
the sigma function (see 
$\S$20.53 in \cite{ww}). If we further assume that we have integer initial 
data satisfying (\ref{edsinits}), then the integrality of 
the subsequent terms  
and the divisibility 
property (\ref{divis}) can be proved by induction using the 
relations (\ref{hankel}). For other special properties of integer 
EDS we refer the reader to Ward's papers \cite{ward1, ward2} and to  
the more recent works \cite{eew,ems,recs,shipsey,silverman,swart}.

If we are working with an elliptic curve $E$ over $\C$, then the terms  
$W_n=\si (n\ka )/\si (\ka )^{n^2}$ of 
a generalized EDS are essentially 
just values of the the division polynomials of $E$ 
(see Exercises 24 and 33 in chapter 20 of \cite{ww}, 
and chapter II of \cite{lang}), possibly up 
to  a prefactor of the form $\gamma^{1-n^2}$ (see Lemma 7 in 
\cite{silverman}). Each $W_n$ is an elliptic function of $\ka$ 
that can be written as a polynomial in the variables 
$(\la , \mu )=(\wp (\ka ),\wp '(\ka ))$ with coefficients 
in $\Q [g_2,g_3]$ (cf. Exercise 3.7 in \cite{silver1}, 
or chapter II in \cite{lang}). 
For our purposes, in order to make use of the 
EDS that is the ``companion'' of a Somos 4 sequence, 
we will need to write the $W_n$  as 
polynomials in a slightly different set of variables (Theorem 2.2 below).  

Given any Somos 4 sequence satisfying a bilinear recurrence 
(\ref{bil}) with coefficients $\al$, $\be$, 
there is a natural companion EDS associated 
with it, which can be defined in two different ways 
\cite{swartvdp, hones5}. 

\noindent {\bf Definition 1. Algebraic definition 
of companion EDS:} 
{\it For a Somos 4 sequence, compute the translation invariant 
according to (\ref{tinvt}),  
and then the companion 
EDS $\{ \,W_n \, \}_{n\in\Z}$ 
is the (generalized) EDS  that satisfies the {\it same} fourth-order 
recurrence (\ref{bil}) with initial values specified to be  
\beq 
W_1=1, \quad W_2=-\sqrt{\al}, \quad W_3=-\beta, \quad 
W_4 =\I\sqrt{\al} , 
\label{compeds} 
\eeq 
where 
\beq 
\I =\al^2 +\be T. \label{idef} 
\eeq 
(Note that the coefficients of the 
Somos 4 sequence are then given in terms of the first 
three terms of its companion EDS by (\ref{edsp}).) 
} 

\vspace{.1in}  
\noindent {\bf Definition 2. Analytic definition 
of companion EDS:} 
{\it  
For a Somos 4 sequence, according to 
Theorem 1.2 the terms $A_n$ are given explicitly by 
the formula (\ref{form}) for suitable parameters 
$g_2,g_3,\hat{a},\hat{b},z_0,\ka $, and then 
the terms $W_n$ of the companion EDS are given  in terms of the 
Weierstrass sigma function by the formula (\ref{wardeds}) with the 
same parameters $g_2,g_3,\ka$. } 

\vspace{.1in} 
In his original memoir \cite{ward1}, Morgan Ward showed that 
an EDS admits two equivalent  definitions, one algebraic and the other analytic. 
By construction, the translation 
invariant for the Somos 4 sequence and that for its companion 
EDS have  the same value $T=\wp ''(\ka )$, and, upon 
comparing (\ref{edsp}) with (\ref{parform}), the algebraic and 
analytic definitions of the companion EDS are easily 
seen to be equivalent over $\C$; see also the Remarks  
after equations (1.18) and (2.14) in \cite{hones5}.  
The introduction of the companion EDS is very natural in the light 
of the following result due to one of us with 
van der Poorten \cite{swartvdp}. 

\begin{thm} 
The terms $A_n$ of a Somos 4 sequence satisfy the 
Hankel determinant formulae   
\beq 
\label{hankel2} 
(W_1)^2 A_{n+m} A_{n-m}=   
\left|\begin{array}{cc}   
W_m A_{n+1} & W_{m-1} A_n \\   
W_{m+1}A_n & W_m A_{n-1}   
\end{array} \right|    
\eeq
and 
\beq   
\label{hankel3}  
W_1 W_2 A_{n+m+1}A_{n-m}=  
\left|\begin{array}{cc}  
W_{m+1} A_{n+2} & W_{m-1} A_n \\  
W_{m+2}A_{n+1} & W_m A_{n-1}  
\end{array} \right|
\eeq
for all $m,n\in \mathbb{Z}$, where $W_m$ 
are the terms of the companion EDS. 
\end{thm} 

\begin{prf} 
For a purely algebraic proof of (\ref{hankel2}) and (\ref{hankel3}) 
see \cite{swartvdp}. 
For an analytic proof based on the 
three-term equation for the sigma function see 
Corollary 1.2 and Corollary 1.3 in \cite{hones5}.  
Note that although $W_1=1$ we have included the $W_1$ terms on the 
left hand sides above to illustrate the homogeneity of these 
expressions.   
\end{prf} 

It is clear that a (generalized) EDS is its own companion, and Ward's 
formula (\ref{hankel}) is the special case of (\ref{hankel2}) when 
$A_n=W_n$. Theorem 2.1 is our main tool for proving integrality 
properties of Somos 4 sequences. To begin with we can use it to derive 
a property of EDS which is apparently new.    

\begin{thm} 
{\bf Polynomial representation for elliptic divisibility sequences:}
For a (generalized)   EDS defined by the fourth-order recurrence 
(\ref{eds})  
with 
initial data $W_1=1$, $W_2=-\sqrt{\al}$,  
$W_3=-\be$, $W_4=\I\sqrt{\al}$, the terms in the sequence
satisfy $W_{2n-1} \in \Z [\al^2 , \be , \I ]$ and 
$W_{2n} \in \sqrt{\al}\,\,\Z [\al^2 , \be , \I ]$ for all $n\in \Z$.
\end{thm} 

\begin{prf}Since $W_{-n}=-W_n$ we need only 
consider 
$n>0$, and then 
the proof is by induction on $n$. Clearly the result is true 
for the initial data $W_1,W_2,W_3,W_4$, 
and the next terms are 
$$ W_5=-\al^2\I +\be^3, 
\qquad 
W_6=-\sqrt{\al}(\I^2 +\al^2 \I -\be^3)\be ,  
$$ 
so the appearance of $\al^2$ is evident. Now setting 
$n=m+1$ and $n=m+2$ in  Ward's formula 
(\ref{hankel}), and putting in the initial values 
$W_1$ and $W_2$,  yields  
\beq 
\label{h1}  
W_{2m+1}=(W_m)^3W_{m+2}-(W_{m+1})^3W_{m-1} 
\eeq 
and 
\beq 
\label{h2} 
W_{2m+2}=\frac{(W_{m+2})^2W_{m+1}W_{m-1} 
-(W_m)^2W_{m+3}W_{m+1}}{\sqrt{\al}} 
\eeq  
respectively, which correspond to well known identities  
for division polynomials -  see e.g. Exercise 3.7 in \cite{silver1}. 
Each term on the right hand side of 
(\ref{h1}) is a product of either four even index terms or four odd index 
terms, and hence 
by the inductive hypothesis, for $m\geq 3$, 
$W_{2m+1}\in\Z [\al^2 , \be , \I ]$. Similarly 
for $W_{2m+2}$, 
both terms in the numerator on the right hand side of
(\ref{h2}) consist of a product of two odd index and two even index 
terms, and so lie in $\al \, \Z [\al^2 , \be , \I ]$ by the inductive 
hypothesis; dividing out by $\sqrt{\al}$ in the denominator gives the 
required result.     
\end{prf} 

\noindent {\bf Remarks.} Taking the Weierstrass model 
\beq \label{wab} 
Y^2=X^3+AX+B, 
\eeq 
which is equivalent to (\ref{wcurve}) via 
$x=X$, $y=2Y$, $g_2=-4A$, $g_3=-4B$, 
it is known that the division polynomials corresponding to the 
multiples $[n](X,Y)$ are elements of $\Z [A,B,X,Y^2]$ for $n$ odd and 
of $2Y\,\Z [A,B,X,Y^2]$ for $n$ even. Hence using 
the equation (\ref{wab}) to eliminate $Y^2$, and dividing by $2Y$ where 
necessary,  gives 
polynomials in  $\Z [A,B,X]$ 
(see Exercise 3.7 in \cite{silver1}, or \cite{lang}, chapter II). Thus we see that the result of 
Theorem 2.2 corresponds to a different choice of basis for 
the division polynomials: after suitable rescaling by $2Y$ for 
$n$ even,  they are polynomials in the variables 
$$ 
\al^2=16Y^4, \qquad \be = A^2-6X^2A-12XB -3X^4,
$$ 
$$ 
\I = 
2A^3+10X^2A^2-(10X^4-8XB)A -2X^6-40X^3B+16B^2
$$ 
with integer coefficients. 
As far as we are aware, this is a new result about
division polynomials.
The quantities $\be ,\I$  and (via (\ref{wab})) $\al^2$ are 
themselves elements of $\Z [A,B,X]$.  
  
Theorems 2.1 and 2.2 together imply a stronger version of the 
Laurent property for Somos 4 sequences, 
which lead us to establish a fairly weak set of criteria for 
integrality in the next section. 

\section{Proof of the strong Laurent property} 

\setcounter{equation}{0}

With the above results at hand, we are now able to 
state our main result concerning a strong version of the Laurent property for 
Somos 4 sequences.

\begin{thm}\label{strongs4}
{\bf The strong Laurent property for Somos 4 sequences:}
Consider a Somos 4 sequence defined by a fourth-order recurrence
(\ref{bil}) with coefficients $\al$, $\be$,  
initial data $A_1, A_2, A_3, A_4$, 
and translation invariant $T$ as in (\ref{tinvt2}). 
The subsequent terms of the sequence
are elements of the ring $\mathcal{R}$ 
of polynomials 
generated by these coefficients and initial data, as well 
as $A_1^{-1}$ and 
the quantity $\I =\al^2 +\be T$. In other words,   
$$A_n \in
\mathcal{R}=\Z [\al , \be , \I , A_1^{\pm 1},
A_2,  A_3,  A_4]$$ 
for all $n\geq 1$.
\end{thm}

\begin{prf} 
Setting $n=m+1$ in each of the Hankel determinant formulae 
(\ref{hankel2}) and (\ref{hankel3}) and substituting $W_1=1$, 
$W_2=-\sqrt{\al}$ leads to 
\beq
\label{d1}
A_{2m+1}=\frac{(W_m)^2A_mA_{m+2}-W_{m+1}W_{m-1}A_{m+1}^2}{A_1}  
\eeq
and
\beq
\label{d2}
A_{2m+2}=\frac{W_{m+2}W_{m-1}A_{m+1}A_{m+2}  
-W_mW_{m+1}A_mA_{m+3}}{\sqrt{\al}A_1}. 
\eeq
The proof that $A_n\in\mathcal{R}$ for $n\geq 1$  
then proceeds by induction, using the identities (\ref{d1}) 
and (\ref{d2}) for odd/even $n$ respectively. For odd $n=2m+1$ 
it is clear from Theorem 2.2 that, on the 
right hand side of (\ref{d1}) each of the 
terms $W_m^2$ and   $W_{m+1}W_{m-1}\in\Z [\al , \be ,\I ]$, and by the 
inductive hypothesis the other terms in (\ref{d1}) are in 
$\mathcal{R}$. Similarly, 
for even $n=2m+2$, on the right hand side 
of (\ref{d2}) it is clear that  both $W_{m+2}W_{m-1}/\sqrt{\al}$ and    
$W_mW_{m+1}/\sqrt{\al}\in\Z [\al , \be ,\I ]$, and the result follows. 
\end{prf} 

\noindent {\bf Remarks.} From (\ref{tinvt2}) it is clear that 
$T\in \Z [\,\al , \be , A_1^{\pm 1}, 
A_2^{\pm 1}, A_3^{\pm 1}, A_4^{\pm 1} \, ]$, 
so 
$\I \in \Z [\,\al , \be , A_1^{\pm 1},  
A_2^{\pm 1}, A_3^{\pm 1}, A_4^{\pm 1} \, ]$.  
Similarly, applying 
Theorem 3.1 to the reversed sequence 
$A_n^*=A_{5-n}$ starting from 
$A_1^*=A_4$ gives 
$A_n \in
\Z [\al , \be , \I , A_1,
A_2,  A_3,  A_4^{\pm 1}]$ 
for all $n\leq 4$.
Hence it follows that for all $n\in\Z$, the terms 
$A_n \in \Z [\al , \be , A_1^{\pm 1},
A_2^{\pm 1},  A_3^{\pm 1},  A_4^{\pm 1}]$, 
so Theorem 1.1 is a consequence of Theorem 3.1. However, 
note that as it stands, Theorem 3.1 is unidirectional (it 
only applies for $n\geq 1$), whereas Theorem 1.1 applies 
to both positive and negative $n$.  

\begin{cor} {\bf Integrality criteria for 
Somos 4 sequences:} \label{strongcors4}
Suppose that a Somos 4 sequence is defined by the 
recurrence (\ref{bil}) for coefficients 
$\al ,  \be\in\Z$, and initial 
values $A_1=\pm 1$ with non-zero $A_2,A_3,A_4\in\Z$. If $\be T$ 
is an integer, where $T$ is the translation 
invariant given by (\ref{tinvt2}), 
then $A_n\in\Z$ for all $n\geq 1$. If it further holds that   
the backwards iterates $A_0,A_{-1},A_{-2}\in\Z$, then 
$A_n\in\Z$ for all $n\in\Z$. 
\end{cor}

\noindent {\bf Proof of Corollary \ref{strongcors4}.} For $n\geq 1$ the integrality 
is obvious, since $\I =\al^2+\be T$ and 
$\al$, $\be\in\Z$. 
If the three backward iterates adjacent to $A_1$ are 
also integers, then to get integrality for $n<1$  
it suffices to apply Theorem 3.1 to the 
reversed sequence $B_n=A_{2-n}$, and then $B_n\in\Z$ for 
$n\geq 1$ and the result follows. $\Box$

\begin{cor}\label{quadioph} 
Given $\al ,\be \in \Z$ and $T\in\Q$ with $\be T\in\Z$, 
suppose that the quartic equation 
\beq\label{s4dio} 
s^2 v^2+\al \, (su^3+ t^3v) 
+ \be \,  t^2 u^2
=T\,stuv 
\eeq 
has a solution of the form $(s,t,u,v)=(A_1,A_2,A_3,A_4)$, with $A_1=\pm 1$ 
and non-zero  
$A_2,A_3,A_4\in\Z$. Then provided that the orbit of this set of initial data 
under (\ref{bil})
is non-periodic, 
it produces infinitely many 
integer solutions of the Diophantine equation (\ref{s4dio}).    
\end{cor} 

\noindent {\bf Remarks.} 
Note that the equation (\ref{s4dio}) is homogeneous in 
$(s,t,u,v)$, so it is really the solutions with $\gcd (s,t,u,v)=1$ that are 
the interesting ones. 
For the sequence (\ref{s4seq}) it is easy to show by induction that 
$\gcd  (A_n,A_{n+1},A_{n+2},A_{n+3}) =1$ for all $n$;  
but in general for sequences where this holds for $n=1$ 
then it need not be so for all $n$, particularly when 
$\gcd (\al ,\be )\neq 1$; see Example \ref{exa1} 
for instance. Nevertheless, in such cases it can still be often be 
checked that one does get infinitely many distinct  solutions 
satisfying $\gcd (s,t,u,v)=1$.

\noindent {\bf Proof of Corollary \ref{quadioph}.} From Theorem \ref{swahonresults} 
we see that 
the equation (\ref{s4dio}) is just a rewriting of (\ref{tinvt}). Since $T$ is a conserved 
quantity for (\ref{bil}), for all $n$ the quadruple 
$(s,t,u,v)=(A_n,A_{n+1},A_{n+2},A_{n+3})$ lies on 
the quartic threefold defined by (\ref{s4dio}), and by Corollary \ref{strongcors4} this is 
a quadruple of integers for $n\geq 1$.   
As long as the orbit of the initial quadruple $(A_1,A_2,A_3,A_4)$ is not periodic 
(which would correspond to $P\in E$ being a torsion point), these quadruples are all 
distinct.  
$\Box$

\begin{exa} \label{exa1}
Consider the 
sequence defined by 
$$ 
A_{n+4}A_{n} = 11^3  A_{n+3}A_{n+1} + 90\cdot 11^3 A_{n+2}^2  
$$ 
with initial values
$A_1=1, A_2=3, A_3=11^2, A_4=11^3\cdot 7\cdot 19$; this is 
just the example previously mentioned  with (\ref{1stex}). This 
sequence extends in both directions, thus: 
\beq \label{sbag} 
\ldots ,2498287 
,1221,7,1,3,121,177023,2460698229,
\ldots \eeq 
In this case we have  integer coefficients 
$\al =1331$, $\be =119790$, and with $n=2$ in the formula 
(\ref{tinvt}) from Theorem 
\ref{swahonresults} we calculate 
$T=869$, so certainly $\be T\in \Z$. Also we have $A_1=1$, 
and from (\ref{sbag}) we see that on each side of this, 
the three adjacent iterates ($A_2,A_3,A_4$ and $A_0,A_{-1},A_{-2}$ 
respectively) are all integers, so all the   conditions of Corollary 
\ref{strongcors4} are met and this sequence 
consists entirely of integers. It is easy to see that the terms of the sequence 
grow rapidly,  such that $\log A_n \sim C n^2$ with $C\approx 1.5$. 
So the sequence cannot be periodic, and hence (by Corollary \ref{quadioph})  
it yields infinitely many solutions of the Diophantine equation
$$ 
s^2 v^2+1331(su^3+ t^3v ) +  119790 t^2 u^2=869stuv.  
$$ 
\end{exa} 

The preceding example admits a rather broad generalization, as follows. 

\begin{exa}\label{exabcde} 
Pick three integers $a,d,e$, and 
choose an ordered pair of integer factors $b,c$ such that 
\beq \label{quarticvariety}
a^3 \, d +e^2= b \, c.
\eeq 
Now take the Somos 4 recurrence of the form  
\beq \label{ades4}
A_{n+4} A_{n}= e^3 \, A_{n+3} A_{n+1} + a d e^3 \, (A_{n+2})^2, 
\eeq 
with the initial conditions $A_1=1$, $A_2=a$, $A_3=e^2$, $A_4=ce^3$. 
This extends in both directions as 
\beq \label{abcde} 
\ldots , e^3( b^2 d + e(b+d) ) , ae(b+d), b, 1, a, e^2, ce^3, ae^6(c+de), \ldots 
\eeq 
so that the three adjacent iterates both to the left and to the right 
of $A_1$ are integers. Thus we have $\al = e^3$, $\be = ade^3\in\Z$ and 
from (\ref{tinvt}) we calculate $\be T=de^4 (a^3+be+c) \in \Z$, which means 
that  Corollary 
\ref{strongcors4} applies and $A_n\in\Z$ for all $n\in\Z$. 
\end{exa} 

\noindent {\bf Remarks.} The preceding example can be interpreted in several 
different ways. Working over $\C$, the equation (\ref{quarticvariety}) 
defines an affine variety $V$ in $\C^5$, and then by Theorem \ref{swahonresults} 
there is a family of elliptic curves fibred over $V$, with the sequence 
(\ref{abcde}) corresponding to a sequence of points on each fibre. 
Note also that by setting $X=-ad$, $Y=de$, $D=bcd^2$ we get 
\beq \label{dcurve}
Y^2=X^3+D, 
\eeq 
so that the variety $V$ itself admits a fibration with all generic 
fibres ($bcd\neq 0$) being isomorphic over $\C$ to the same elliptic curve 
(since all the curves (\ref{dcurve}) are isomorphic for $D\neq 0$). 
Alternatively, we can consider the function field generated by the 
variables $a,b,c,d,e$ subject to the relation (\ref{quarticvariety}). In that case, 
Theorem \ref{strongs4} implies that the recurrence (\ref{ades4}) generates 
a sequence of polynomials in the ring $\Z [a,b,c,d,e] /{}\sim$, where 
$\sim$ is the corresponding equivalence relation. Note also 
that Example \ref{exa1} is the particular case $a=3$, $b=7$, $c=133$, 
$d=30$, $e=11$, while the original Somos (4) sequence given by  
(\ref{s4seq}) is the case $a=1$, $b=2$, $c=1$, $d=1$, $e=1$. 
Generically (provided that $b\neq c$ and $e\neq 1$) a different 
ordering of the factors 
$b,c$ in (\ref{quarticvariety}) not only gives a different 
sequence (\ref{abcde}) but also corresponds to a sequence of 
points on a different (non-isomorphic) elliptic curve.

We should point out that the sufficient criteria for
integrality in
Corollary \ref{strongcors4} are not necessary conditions. 
To see this, it is enough to 
consider the fact that
Somos 4 sequences are invariant under the two-parameter
abelian group of gauge transformations defined by
\beq \label{igauge}
A_n\longrightarrow \tilde{A}_n=\tilde{a}\,\tilde{b}^n\,A_n, \qquad
\tilde{a},\tilde{b} \in \C^*,
\eeq
in the sense that $\tilde{A}_n$ satisfies the same bilinear
recurrence (\ref{bil}) as does $A_n$. Furthermore
applying (\ref{igauge}) leaves the
the translation invariant $T$ the same,
and hence the curve $E$ in Theorem \ref{swahonresults}
is preserved by the action of this gauge group. Moreover,
if $A_n\in\Z$ for all $n$ and the transformation (\ref{igauge})
is applied with $\tilde{a}$, $\tilde{b}\in\Z^*$,
then integrality is preserved for $n\geq 1$, and if we fix
$\tilde{b}=1$ then
$\tilde{A}_n\in\Z$ for all $n$, and $|\tilde{A}_n|\geq 2$ for
$|\tilde{a}|\geq 2$. 

As well as the gauge transformations (\ref{igauge}), Somos 4 sequences
exhibit an additional 
{\it covariance} property: under
the action of
the group of transformations given by  
\beq \label{cgauge}
A_n\longrightarrow \hat{A}_n=\hat{c}^{n^2}\,A_n, \qquad
\hat{c} \in \C^*,
\eeq
the terms $\hat{A}_n$
of the transformed sequence satisfy a Somos 4 recurrence, 
which is of the same form
(\ref{bil}) but with the coefficients
and translation invariant rescaled according to
\beq \label{covab}
\hat{\al} = \hat{c}^6\al , \qquad \hat{\be}= \hat{c}^8\be,
\qquad \hat{T}=\hat{c}^4T,
\eeq
so that the $j$-invariant (\ref{jinvt}) is preserved.
For each $\hat{c}$, the transformation (\ref{cgauge}) just corresponds
to homothety acting on the Jacobian: $z\in \mathrm{Jac}(E)$ is mapped to
$\hat{z}=\hat{c}^{-1}z\in \mathrm{Jac}(\hat{E})$. It is clear that
applying the transformation (\ref{cgauge}) with
$\hat{c}\in\Z^*$ preserves integrality of
Somos 4 sequences,
as does the transformation
$$
A_n\longrightarrow A_n^\dagger =(A_k)^{(n-k)^2-1}\, A_n
$$
which inserts a $1$ at the $k$th term of the sequence.
All transformations of this type, that include elements of the
group (\ref{cgauge}) with $\hat{c}\in\Z^*$,  have the effect
of increasing the discriminant of
the associated curve $E$ when $|\hat{c}|>1$. 

Corollary \ref{strongcors4} only applies when one of the 
terms in a given sequence is $\pm 1$ and sufficiently many 
adjacent terms are integers. If this is not the case, then one might 
try to engineer it to be so 
by applying a combination of the transformations (\ref{igauge}) and 
(\ref{cgauge}), in order to minimize the size of the discriminant as much as 
possible while preserving the requirement that 
$\al$, $\be$, $\I \in \Z$ as well the the integrality of a given set of adjacent terms. 
However, without the requirement that one of the terms of the 
sequence should be $\pm 1$, we 
can present a different set of sufficient criteria for integrality.   

\begin{thm}  
If $\alpha, \beta\in\Z$ and eight successive terms 
$A_{-2},A_{-1} \ldots, A_4,  
A_5$ are integers, and if 
$$\gcd(\alpha, \beta) = \gcd(A_1,A_2) =  
\gcd(\alpha,A_0,A_2) = \gcd(\alpha, A_1,A_3) = 1,$$ 
then the terms $A_n$  
of a Somos 4 sequence are integers for all $n\in\Z$.  
\end{thm}  
  
\begin{prf}  
Suppose that a prime $p$ appears in the denominator of some term,  
and suppose that $p\not |A_1$.  Then by Theorem 3.1, $p$  
divides the denominator of $\I$, and hence it must divide the 
denominator of $T\in\Q$. 
Setting $n=1$ and $m=3$ 
in the result (\ref{hankel2}) 
from Theorem 2.1,  we have   
$$ 
A_{-2} A_4  
-\be^2 A_0 A_2 = \I \alpha (A_1)^2\in\Z      
$$ 
by the integrality assumptions, whence $p| \al$; 
and then $p\not| \be$ because $\gcd (\al , \be ) =1$. 
Thus, making use of the recurrence 
(\ref{bil})  gives 
$$ 
A_{-1} A_3 = \alpha A_0 A_2  
+ \beta (A_1)^2 \equiv \beta (A_1)^2 \not\equiv 0 \bmod p 
$$ 
But  
from the expression (\ref{tinvt3}) with $n=0$ and $n=2$ we have that 
$p|A_{-1}A_0A_1$ and $p|A_1A_2A_3$ respectively. 
It follows 
that $p$ divides $\gcd(A_0, A_2)$ and hence $\gcd(\alpha, A_0,  
A_2)$. But $\gcd(\alpha, A_0,A_2)=1$, contradicting the initial 
assumption that $p\not|A_1$.   

Now by repeating the above argument with  
each $A_n$ replaced by $A_{n+1}$, it is clear that if  
$p$ is a prime which appears in the denominator of some term, and if  
$p$ is coprime to $A_2$, then $p$ divides $\gcd(\alpha, A_1,  
A_3)$. But then, since $\gcd(\alpha, A_1,A_3)=1$,  
we must have also that $p|A_2$. Thus we require that 
$p|\gcd(A_1,A_2)$ and the result follows by contradiction.  
\end{prf}  


It has been known for some time that Somos 5 sequences, 
defined by bilinear recurrences 
of the form   
\beq 
\tau_{n+3}\tau_{n-2}=\tal \,\tau_{n+2}\tau_{n-1}+\tbe \,\tau_{n+1}\tau_n,  
\label{s5} 
\eeq 
are also related to sequences of points on elliptic curves 
\cite{zagier, vdp, swartvdp} (see also the 
unpublished results of Elkies quoted in 
\cite{heron}). The complete solution of the initial value problem for the 
general    
Somos 5 sequence, in terms of the Weierstrass sigma function, 
was first presented in \cite{hones5}. It turns out that Somos 5 sequences 
have two algebraically independent analogues of the translation 
invariant, but for our purposes here we will only need 
the quantity $\tJ$ from \cite{hones5} given in terms of five initial 
data by 
\beq \label{tj} 
\tJ= 
\frac{\tau_{4}\tau_{1}}{\tau_{3}\tau_{2}} + 
\frac{\tau_{5}\tau_{2}}{\tau_{4}\tau_{3}} +  
\tal\left( 
\frac{\tau_{3}\tau_{2}}{\tau_{4}\tau_{1}} 
+ 
\frac{\tau_{4}\tau_{3}}{\tau_{5}\tau_{2}} \right) 
+\tbe\frac{(\tau_{3})^2} {\tau_{5}\tau_{1}}  
. 
\eeq 
Rather than summarizing  all the results of \cite{hones5}, we will 
just state our  main theorem concerning the Laurent 
property and refer to some of these results as
needed in the proof.  

\begin{thm} \label{strongs5}
{\bf The strong Laurent property for Somos 5 sequences:}
For a Somos 5 sequence defined by a fifth-order recurrence
(\ref{s5}) with coefficients $\tal$, $\tbe$ and  
initial data $\tau_1, \tau_2, \tau_3, \tau_4, \tau_5$,
the quantity 
\beq \label{5jinvt} 
\tJ=\frac{\tau_{n+1}\tau_{n-2}}{\tau_{n}\tau_{n-1}} + 
\frac{\tau_{n+2}\tau_{n-1}}{\tau_{n+1}\tau_{n}} + 
\tal\left( 
\frac{\tau_{n}\tau_{n-1}}{\tau_{n+1}\tau_{n-2}} 
+ 
\frac{\tau_{n+1}\tau_{n}}{\tau_{n+2}\tau_{n-1}} \right) 
+\tbe\, \frac{(\tau_{n})^2}{\tau_{n+2}\tau_{n-2}}
\eeq 
is independent of $n$. 
The subsequent terms of the sequence
are elements of the ring of polynomials
generated by these coefficients and initial data, as well
as $\tau_1^{-1}$, $\tau_2^{-1}$ and
the quantity $\tI =\tbe +\tal \tJ$. 
In other words,
$$\tau_n \in
\Z [\tal , \tbe , \tI  , \tau_1^{\pm 1},
\tau_2^{\pm 1},  \tau_3,  \tau_4, \tau_5]$$ for all $n\geq 1$.
\end{thm}

\begin{prf} 
The fact that $\tJ$ given by (\ref{5jinvt}) 
is independent of $n$ is part of Theorem 2.5 in 
\cite{hones5}, where it is rewritten in terms 
of the quantity $$h_n=\tau_{n+2}\tau_{n-1}/(\tau_{n+1}\tau_n).$$ 
According to Corollary 2.10 in \cite{hones5} 
(or the comments in Section 7 of \cite{swartvdp}), each 
Somos 5 sequence also has a companion EDS associated with it, with 
terms denoted $a_n$ there, such that the 
analogue of (\ref{hankel3}), given by 
\beq 
\label{hankel4}
a_1 a_2 \tau_{n+m+1}\tau_{n-m}=
\left|\begin{array}{cc}
a_{m+1} \tau_{n+2} & a_{m-1} \tau_n \\
a_{m+2}\tau_{n+1} & a_m \tau_{n-1}
\end{array} \right|
\eeq 
holds for all $m$, $n\in\Z$. 
Using further results of that Theorem 2.5, we 
can define this companion EDS by the Somos 4 recurrence 
$$ 
a_{n+4}a_n=\tmu^2\, a_{n+3}\,a_{n+1}-\tal \,(a_{n+2})^2 
$$ 
with initial data 
$ 
a_1=1$,  $a_2=-\tmu$, $a_3=\tal$,  
$a_4=\tmu\tbe$, 
where 
\beq \label{mudef} 
 \tmu = (\tbe +\tal\tJ )^{1/4}. 
\eeq  
It follows from Theorem 2.2, making the 
replacements $\sqrt{\al}\to \tmu$, $\be\to -\tal$, 
$\I\to\tbe$, that $a_{2n+1}$ and $a_{2n}/\tmu$ are 
both elements of $\Z [\tmu^4,\tal ,\tbe ]$, so by   
(\ref{mudef}), noting that $\tI =\tmu^4$  we have 
$a_{2n+1}\in\Z [\tal ,\tbe ,\tI ]$ and 
$a_{2n}\in \tmu\,\Z[\tal ,\tbe ,\tI ]$ for all $n$. Setting 
$n=m+1$ and $n=m+2$ in (\ref{hankel4}) yields two identities 
for $\tau_n$ which for even/odd index $n$ have  
respectively $\tau_1$ and  $\tau_2$ 
in the denominator. Thereafter the proof proceeds 
similarly 
to that of Theorem \ref{strongs4}.   
\end{prf} 
 
\noindent {\bf Remarks.} Similarly to 
the Remarks after Corollary \ref{strongcors4}, the Laurent property for Somos 5 
sequences, as proved in \cite{fz}, is a consequence of Theorem \ref{strongs5}. 
Also, the latter has two immediate corollaries which are the respective analogues of 
\ref{strongcors4} and  \ref{quadioph}.

\begin{cor} {\bf Integrality criteria for 
Somos 5 sequences:} \label{strongcors5} 
Suppose that a Somos 5 sequence is defined by the
recurrence (\ref{s5}) for coefficients
$\tal ,  \tbe\in\Z$, and initial
values $\tau_1=\pm 1$, $\tau_2=\pm 1$, with non-zero  
$\tau_3,\tau_4,\tau_5\in\Z$. 
If the quantity $\tal\tJ$ is an integer, with $\tJ$  
given by (\ref{tj}), 
then $\tau_n\in\Z$ for all $n\geq 1$. If it further holds that
the backward iterates $\tau_0,\tau_{-1},\tau_{-2}\in\Z$, then
$\tau_n\in\Z$ for all $n\in\Z$.
\end{cor} 

\begin{cor}\label{quindioph} 
Given $\tal ,\tbe \in \Z$ and $\tJ\in\Q$ with $\tal\tJ\in\Z$, 
suppose that the quintic equation 
\beq\label{s5dio} 
(\tis\tiw + \tal \, \tiu^2)(\tis\tiv^2 + \tit^2\tiw )  
+ \tbe \,  \tit\tiu^3\tiv  
=\tJ\,\tis\tit\tiu\tiv\tiw  
\eeq 
has an integer
solution of the form $(\tis ,\tit ,\tiu ,\tiv ,\tiw )=(\tau_1,\tau_2,\tau_3,\tau_4,\tau_5)$, with $|\tau_1|= 1$, 
$|\tau_2|= 1$ 
and non-zero  
$\tau_3,\tau_4,\tau_5\in\Z$. Then provided that the orbit of this set of initial data 
under (\ref{s5})
is non-periodic,  
it produces infinitely many 
integer solutions of the Diophantine equation (\ref{s5dio}).    
\end{cor}

\begin{exa} \label{sbagliato} 
Take $\tal = 14641=11^4$, $\tbe =  1771561=11^6$ in (\ref{s5}) with 
initial data $847,8,1,1,33$. From (\ref{5jinvt}) we calculate 
$\tJ = 627$, hence $\tal\tJ\in \Z$, and the sequence extends on either side of these 
values as 
\beq\label{sbagorig}
\ldots ,  805255, 847,8,1,1,33, 6655, 19487171,\ldots, 
\eeq 
so all the conditions of Corollary \ref{strongcors5} hold. 
Thus $\tau_n\in\Z$ for all $n\in\Z$. By Corollary \ref{quindioph}, 
for these values of $\tal$, $\tbe$, $\tJ$ the sequence  
(\ref{sbagorig}) gives  
infinitely many quintuples of integer solutions of   
the corresponding Diophantine equation, of the form (\ref{s5dio}).   
\end{exa} 

\noindent {\bf Remarks.} It is shown in \cite{hones5} (Proposition 2.8) 
that the subsequences of a Somos 5 sequence consisting of the 
even/odd index terms respectively satisfy the same Somos 4 recurrence 
relation. If we index the terms of the preceding example such 
that $\tau_1 =8$, $\tau_2 =1$, $\tau_3=1$, $\tau_4=33$ etc. 
then we find that 
the even index subsequence of (\ref{sbagorig}) is essentially the same as 
the sequence (\ref{sbag}) in Example \ref{exa1}, up to 
applying transformations of the form 
(\ref{igauge}) and 
(\ref{cgauge}). To be precise, we have 
$\tau_{2n}=11^{\mathcal{Q}(n)}A_n$ where $\mathcal{Q}(n)=(n-1)(3n-4)/2$. 
Indeed, both sequences correspond to a sequence of points on 
the same elliptic curve with $j$-invariant (from (\ref{jinvt})) 
$$ 
j=\frac{5^3 \cdot 23^6 \cdot 1013^3 }{ 2^4 \cdot 11^7 \cdot 17^2 \cdot 37 \cdot 1069} .
$$
It is clear from the denominator that the curve has bad reduction 
$\bmod {}11$,
and all terms of the sequence reduce to $0\bmod 11$ 
apart from the central terms 
$8,1,1$.

Note that so far all of our criteria for integrality have been stated 
for Somos recurrences with integer coefficients. However, even 
when $\al$, $\be \in\Z$ the associated curve (\ref{wcurve}) 
must be  defined over $\Q$, or more generally over $\Q (\sqrt {\hat{\al}})$ 
where $ \hat{\al}$ is the square-free part of $\al$, in order 
to consider the sequence of points (\ref{pts}). Hence it is natural to 
define Somos sequences over $\Q$, which leads to the question of 
whether integer sequences can arise 
for non-integer $\al$ or $\be\in\Q$. The following 
example 
illustrates that this is indeed the case.


\begin{exa} \label{somosexa}
Consider the Somos 4 sequence defined by
(\ref{bil}) with
$\al =-1/2$, $\be = 1$ and    
$A_1=1, A_2=-2 , A_3=2, A_3=1$. 
This sequence has terms 
\beq \label{somcha} 
\ldots ,2, 10,-4, 2,3,1,-2 ,2,1,5,2, 12, -26, 34, 236, 352, -1912, 
\ldots , \eeq 
and 
consists entirely of integers. 
\end{exa}

The sequence (\ref{somcha}) was presented in the Robbins (bilinear) forum 
\cite{propp}
by Michael Somos, who found that  
it should correspond to 
a sequence of points on an elliptic curve $E$ with  $j$-invariant  
\beq 
\label{scj} 
j = \frac{3^3 \cdot   19051^3}{  2^{17} \cdot  1721 };  
\eeq 
this was confirmed by Elkies 
\footnote{See {\tt http://www.math.wisc.edu/\~{}propp/somos/elliptic}}.     
Although Theorem \ref{strongs4} does not apply directly in this 
case, we show in the next section how 
a slight modification of the methods used previously  
provides a proof of integrality for a one-parameter family of Somos 4 sequences 
that generalizes Example \ref{somosexa}.

\section{One-parameter family of integer sequences}

\setcounter{equation}{0}

In this section we consider a one-parameter family of sequences 
\beq 
\label{nseq} 
\ldots, 
1,-N,N,1,1+N^2,N,N^3+2N,-N^4-2N^2-2,N^4+4N^2+2, 
\ldots\eeq 
which is defined by the Somos 4 recurrence (\ref{bil}) 
with coefficients and initial data given by 
\beq \label{ndata} 
\al =-1/N, \quad \be =1, \quad A_1=A_4=1, \quad A_2=-A_3=-N. 
\eeq  
When $N=2$ this is the sequence (\ref{somcha}) in Example 
\ref{somosexa} above. 
The integrality 
properties of the sequence (\ref{nseq}) 
do not follow from Theorem \ref{strongs4},  
because the parameter $\al$ is non-integer for integer $N\neq\pm 1$. 
Nevertheless, 
we 
are able to prove that 
the sequence (\ref{nseq}) consists entirely of polynomials in 
$N$ with integer coefficients, and so the integrality of this sequence 
when $N\in\Z$ follows immediately. From (\ref{tinvt}) 
we 
compute the translation invariant     
$ 
T=-N^2-1/N^2$, and hence from 
(\ref{idef}) 
we have 
$\I =-N^2$. Thus we see that $\I\in\Z$ for $N\in\Z$, which leads to 
special properties for the companion EDS associated with (\ref{nseq}).

\begin{thm}
For the companion EDS associated with (\ref{nseq}), 
defined by the fourth order recurrence
(\ref{eds})
with 
initial data $W_1=1$, $W_2=-iN^{-1/2}$,
$W_3=-1$, $W_4=-iN^{3/2}$, the terms in the sequence
satisfy $W_{2n-1} \in \Z [N^4]$, $W_{4n} \in iN^{3/2}\,\,\Z [N^4]$ and
$W_{4n+2} \in iN^{-1/2}\,\,\Z [N^4]$ for all $n\in \Z$.
\end{thm}

\begin{prf} This follows the same pattern as the proof of 
Theorem 2.2, using the two identities (\ref{h1}) and 
(\ref{h2}), except that in (\ref{h2}) it is necessary to consider 
the cases of odd/even $m$ separately.  
\end{prf} 

\begin{thm} 
Consider the Somos 4 sequence defined by the recurrence (\ref{bil}) 
with coefficients $\al =-1/N$, $\be =1$ and 
initial data 
$A_1=1$, $A_2=-N$, $A_3=N$, $A_4=1$. 
This sequence
satisfies 
$A_n\in\Z [N]$ 
for all 
$n\in\Z$. 
Moreover, $A_n\in \Z[N^2]$ whenever $n\equiv 0$ or $1$ mod $4$, 
while $A_n\in N\,\Z[N^2]$ whenever $n\equiv 2$ or $3$ mod $4$.   
\end{thm} 

\begin{prf} 
This proceeds by induction for $n\geq 1$, analogously to the proof of 
Theorem 3.1, using the two identities (\ref{d1}) and 
(\ref{d2}), except that each of these formulae must be considered  
separately for odd/even $m$ to get the different properties 
of $A_n$ as $n$ varies mod 4. We omit further details, 
except to mention that to cover $n<1$ it suffices to consider 
the reversed sequence $\hat{A}_n=A_{5-n}$ and then proceed as before.  
\end{prf} 

\begin{cor} 
For all $N\in \Z$, the sequence 
(\ref{nseq}) defined by the Somos 4 recurrence (\ref{bil}) with 
coefficients and initial data given by (\ref{ndata}) is an 
integer sequence.     
\end{cor} 

Thus we see that the sequence
(\ref{somcha}) is the particular case $N=2$ of 
the family of integer sequences given by (\ref{nseq}) 
with $N\in\Z$. It is then 
interesting to consider the family of elliptic curves 
corresponding to these sequences.    
By Theorem 1.2, using the formulae (\ref{g2}), 
(\ref{g3}) 
and (\ref{pts}) from section 2,  
for each $N\neq 0$ the sequence (\ref{nseq}) is associated with the sequence 
of points  
$$ 
\left(\,  
-\frac{(N^4-1)^2}{12N^3}  
-\frac{N}{(N^2-1)^2} \, ,\, 
\frac{2iN\sqrt{N}}{(N^2-1)^3}\, \right) + [n]  
\left(\, -\frac{(N^4-1)^2}{12N^3}  , \frac{i}{ \sqrt{N}}\, \right) 
$$ 
on the curve  
\beq 
\label{ncurve} 
y^2=4x^3 -g_2(N)\, x-g_3(N) 
\eeq 
with invariants 
$$ 
g_2(N)=\frac{N^{16}-4N^{12}+30N^8+20N^4+1}{12N^6}, 
$$ 
$$ 
g_3(N)=\frac{N^{24}-6N^{20}+51N^{16}-56N^{12}+195N^8+30N^4+1}{216N^9}, 
$$ 
and 
$j$-invariant 
\beq \label{jn} 
j(N)=\frac{(N^{16}-4N^{12}+30N^8+20N^4+1)^3}{N^{16}(N^{12}-5N^8+39N^4+2)}.   
\eeq 
(The latter reproduces the value (\ref{scj}) when $N=2$.)  

Strictly speaking the recurrence for 
(\ref{nseq}) does not make sense 
when $N=0$, but the sequence (\ref{nseq}) can still be defined 
by setting $N=0$ in each of these polynomials. It can be seen 
that when $N=0$ the values 
of these polynomials are given by 
\beq \label{zvals} 
A_{4m}=(-1)^{m-1}2^{m(m-1)/2} , \,  
A_{4m+1}=2^{m(m-1)/2} , \, A_{4m-1}=A_{4m+2}=0 
\eeq   
for all $m$. The case $N=1$ is also special, because in that 
case the base point $Q=\infty$, and we have  
the sequence 
\beq \label{s4eds} 
1,-1,1,1,2,1,3,-5,7,4,23,29,59,-129,314,\ldots  
\eeq 
with $A_0=0$, which corresponds to 
the multiples $[n]P$ of the point $P=(0,i)$ on 
the curve 
\beq \label{s4curve} 
y^2=4x^3-4x-1, 
\eeq 
for which $j=2^{12}\cdot 3^3/37$. 
In that case, 
the sequence (\ref{nseq}) is almost identical to its companion EDS; 
more precisely, we find that $A_n=(-i)^{n-1}W_n$ for all $n$. Over 
$\C$ (sending $x\to -x$, $y\to iy$), the curve  (\ref{s4curve}) 
is isomorphic to the curve associated with the Somos (4) sequence  
(\ref{s4seq}) as found in \cite{honeblms}, which is not surprising 
when one observes that the terms of (\ref{s4seq}) 
are precisely the odd index terms of the 
sequence (\ref{s4eds}). 
 
The equation (\ref{ncurve}) can be thought 
of as defining a curve $E/\C (N)$. 
Letting $N$ vary this gives an elliptic surface fibred over 
$\mathbb{P}^1$: with the exception of a finite set of values,  
each value of $N$ defines a non-singular  
elliptic curve (see Chapter III in \cite{silver2}). 
The singular fibres correspond to the zeros of the denominator of 
$j(N)$, i.e. $N=0$ and the zeros of $N^{12}-5N^8+39N^4+2$, 
as well as the fibre at infinity.  
If we set $S=N^4$, then it is clear from (\ref{jn}) that 
$j(N)\in \C(S)$. It turns out that 
over $\C (\sqrt{N})$, the rational elliptic surface 
(\ref{ncurve}) is birationally equivalent to the Weierstrass 
model 
\beq 
\label{wmodel} 
\hat{y}^2+8S\,\hat{y}=\hat{x}^3+(S-1)^2\,\hat{x}^2 
-8S(S+1)\,\hat{x} 
\eeq 
which is defined over $\C(S)$. 
In terms of the model (\ref{wmodel}) 
the sequence (\ref{nseq}) corresponds to the sequence of 
points 
$$ 
\hat{Q}+[n]\hat{P}=\left( \, 
\frac{4S}{(\sqrt{S}-1)^2}\, , \, -4S+\frac{8S\sqrt{S}}{(\sqrt{S}-1)^2}\, 
\right) +[n]\Big(\, 0\, , \, 0 \, \Big) 
$$ defined over the extended function field $\C (\sqrt{S})$.

\section{Conclusions} 
\setcounter{equation}{0}

Using the connection between Somos 4 sequences 
and their associated companion EDS, together with 
the identities in Theorem 2.1, we have shown that these sequences 
possess a stronger variant of the Laurent property in \cite{fz}, 
which is expressed by means of the conserved quantity $\I $,  
given by the formula  
(\ref{idef}) in terms of    
$T$,  the translation invariant 
(\ref{tinvt}).  This strong Laurent property has produced a set of sufficient criteria for integrality of 
Somos 4 sequences with integer coefficients, and we 
have given a similar set of criteria for Somos 5 sequences as well. 

The fact that every Somos 4 sequence admits a conserved quantity 
rests on the link with shifts by multiples $[n]P$ 
of a point $P$ on an elliptic curve $E$ \cite{swart}, or equivalently 
with a discrete linear flow on the Jacobian 
$\mathrm{Jac}(E)$. This can also 
be understood from the connection with integrable maps 
\cite{veselov}: 
via the subsitution $f_n=A_{n+1}\,A_{n-1}/(A_n)^2$, 
the bilinear recurrence (\ref{bil}) yields the 
second-order difference equation 
$$ 
f_{n+1}(f_n)^2f_{n-1}=\al f_n +\be, 
$$ 
which has the conserved quantity (\ref{tmap}), 
preserves the symplectic form 
$\omega = (f_{n-1}f_n )^{-1}\, df_{n-1}\wedge df_n$, and is a 
degenerate case of the family of maps studied 
by Quispel, Roberts and Thompson \cite{qrt}. 

There is another aspect of the Laurent property, 
that was stated as a conjecture in \cite{fz}, which 
is that for the octahedron recurrence (also known 
as the discrete Hirota equation \cite{zabrodin},  
an integrable partial difference equation) the 
corresponding Laurent polynomials in the initial data and 
coefficients have all positive coefficients, and furthermore 
all these coefficients are 1. Speyer has given a combinatorial 
proof of this conjecture, using perfect matchings of certain graphs \cite{speyer}, 
and Speyer and Carroll have also used combinatorics to prove the analogous 
result for the four-term cube recurrence (also known as the 
Hirota-Miwa equation). Since the Somos 4 recurrence can be obtained 
as a reduction of the discrete Hirota equation, Speyer notes 
that his result implies that the Laurent polynomials 
in $\Z [\al , \be , A_1^{\pm 1}, 
A_2^{\pm 1},  A_3^{\pm 1},  A_4^{\pm 1}]$ for Somos 4 
have all positive coefficients. 
Whether this result has any interpretation in terms of 
elliptic curves is not clear. 

Although our most general integrality criteria
in Section 3 were formulated for
Somos 4 recurrences with integer coefficients, the 
family
of sequences (\ref{nseq}) illustrates the fact that $\al$ and $\be$
need not be integers for integer sequences
to result.
The sequence of polynomials (\ref{nseq}) 
provides 
further insights into the integer sequences corresponding to 
$N\in\Z$. In particular, consider the case of integer $N>2$, and 
suppose $p>2$ is a prime with $p|N$. In that case, due to the factor 
of $N^{16}$ in the denominator of the $j$-invariant, the curve 
(\ref{ncurve}) - or equivalently (\ref{wmodel}) - has bad reduction 
at the prime $p$. Furthermore, by Theorem 4.2 it is clear that 
$A_{4m+2}\equiv A_{4m+3}\equiv 0$ mod $p$, while from the values 
(\ref{zvals}) of these polynomials at $N=0$ it follows that 
$p\not | A_{4m}$ and $p\not | A_{4m+1}$ for all $m$. 
This means that the distance, or gap length \cite{rob}, 
between successive multiples of  
the prime $p$ appearing in the 
sequence alternates between one and three. This is in contrast 
with Robinson's numerical observations for the Somos (4) sequence 
(\ref{s4seq}), which led him to the conjecture that for this 
particular sequence, the gap length should be {\it constant} 
for any prime $p$ \cite{rob}. The $p$-adic properties 
of EDS have been considered  
recently by Silverman \cite{silverman}, and we expect that 
similar considerations will lead to a better understanding 
of prime divisors in Somos 4 sequences.

Certain higher order Somos 
sequences have also appeared in connection with 
division polynomials \cite{cantor}, higher order integrable maps 
\cite{beh} and divisor sequences for 
hyperelliptic curves \cite{matsupsi, vdp2}. 
Based on na\"{i}ve counting arguments, 
in \cite{honeblms} one of us made the conjecture 
that all Somos sequences should arise from appropriate 
divisor sequences on such curves. 
In the forum \cite{propp}, Elkies has given much stronger 
arguments to the contrary, based on a conjectured 
theta function formula for the iterates. Moreover, the 
Laurent property is failed for Somos 8 recurrences; and for the 
sequence generated by the recurrence 
$$
S_{n+8}S_n=S_{n+7}S_{n+1}+S_{n+6}S_{n+2}+S_{n+5}S_{n+3} 
+( S_{n+4})^2
$$ 
with initial values $S_j = 1$ for $j=1,\ldots ,8$ 
the logarithmic heights of the rational iterates have exponential 
growth, so that $\log h(S_n )\sim K n$ as $n\to\infty$ 
with $K\approx 0.23$; this growth is incompatible with 
a discrete linear flow on an Abelian variety. 
Recently, Halburd has proposed that polynomial growth 
of logarithmic heights should be an integrability criterion for 
birational maps \cite{halburd}.  

One of us  
has found that for Somos 6 (and also Somos 7) recurrences,  
there are two independent conserved 
quantities, 
analogous to the translation invariant (\ref{tinvt}) for Somos 4. 
For these recurrences of sixth and seventh order, it is possible to 
given an explicit formula for the iterates in terms of two-variable theta functions, 
or equivalently in terms of the Kleinian sigma functions for an associated curve 
of genus two. Recently, Kanayama  has derived multiplication 
formulae for genus two sigma functions \cite{kanayama}. One of us has 
pointed out that the statement of Proposition 3 in \cite{kanayama} is 
incorrect, but Kanayama has shown us a corrected version of this result 
\cite{kanayamaprivate}, which we have used 
to  derive  the solutions of Somos 6 and Somos 7 recurrences. 
The full description of the genus two case is the subject of 
further investigation.
    
\noindent {\bf Acknowledgements.} 
AH is grateful to James Propp, Noam Elkies and other members 
of the Robbins bilinear forum \cite{propp} for helpful correspondence
on related matters.

\end{document}